\documentclass[10p]{elsarticle}
\usepackage{amssymb, color}
\usepackage{amsmath}
 \usepackage{pifont}
 \usepackage[utf8]{inputenc}
\usepackage{setspace}
\journal{Journal of Fixed Point Theory and Applications}
\makeatletter
\def\LaTeX{\leavevmode L\raise.42ex
    \hbox{\kern-.3em\size{\sf@size}{0pt}\selectfont A}\kern-.15em\TeX}
\makeatother

\newcommand{\BibTeX}{{\rm B\kern-.05em{\sc
          i\kern-.025emb}\kern-.08em\TeX}}

\makeatletter
\def\@currentlabel{2.1}\label{e:dispaa}
\def\@currentlabel{2.21}\label{e:dispau}
\def\@currentlabel{2.22}\label{e:dispav}
\def\@currentlabel{2.23}\label{e:dispaw}
\def\@currentlabel{2.24}\label{e:dispax}
\def\theequation{\thesection.\@arabic\c@equation}
\makeatother

\renewcommand{\theequation}{\arabic{section}.\arabic{equation}}

\newcommand{\D }{\Delta }

\newcommand{\e }{\varepsilon }

\newcommand{\be}{\begin{equation}}
\newcommand{\ee}{\end{equation}}
\newcommand{\ba}{\begin{align*}}
\newcommand{\ea}{\end{align*}}
\newcommand{\R}{\mathbb R}

\def \D{\Delta}

\newtheorem{thm}{Theorem} [section]
\newtheorem{lem}{Lemma} [section]
\newtheorem{prop}{Proposition} [section]
\newtheorem{cor}{Corollary} [section]
\newtheorem{defn}{Definition} [section]

\newtheorem{re}{Remark}[section]

\renewcommand{\theequation}{\thesection.\arabic{equation}}
\renewcommand{\thesection}{\arabic{section}}
\renewcommand{\theequation}{\thesection.\arabic{equation}}
\let\ssection=\section\renewcommand{\section}{\setcounter{equation}{0}\ssection}

\begin{document}
\begin{frontmatter}

\title{Liouville theorem  on a half-space for biharmonic problem with Dirichlet boundary
condition.}
\author[mf1,mf2]{Foued Mtiri}
\ead{mtirifoued@yahoo.fr}
\author[ah1,ah2]{Abdelbaki Selmi\corref{cor1}}
\ead{Abdelbaki.Selmi@fsb.rnu.tn}
\author[mf3]{ Cherif Zaidi}
\ead{zaidicherif10@gmail.com}
\address[mf1]{Mathematics Department, Faculty of Sciences and Arts, King Khalid University, Muhayil Asir, Saudi Arabia.}
\address[mf2]{ Faculty of Sciences of Tunis, Elmanar University, Tunisia.}
\address[ah1]{ Mathematics Department, Northern Border University, Arar, Saudi Arabia.}
\address[ah2]{Universit\'{e} de Tunis, D\'{e}partement de Math\'{e}matiques, Facult\'{e} des Sciences de Bizerte, Zarzouna, 7021 Bizerte, Tunisia.}
\address[mf3]{Facult\'e des Sciences, D\'epartement de Math\'ematiques, B.P 1171 Sfax 3000, Universit\'e de Sfax, Tunisia .}

\begin{abstract}
\noindent We investigate here the nonlinear elliptic Hénon type equation:
 $$\D^{2} u= |x|^a|u|^{p-1}u \; \,\,\mbox{in}\,\,\,\, \R^{n}_{+},  \quad \quad  u =\frac{\partial u}{\partial x_n} = 0 \quad  \mbox{in}\,\,\,\, \partial \R^{n}_{+},$$
with $p>1$ and $n\geq 2$. In particular, we prove some Liouville type theorems for  stable at infinity solutions. The main methods used are the integral estimates, the Pohozaev-type identity and the monotonicity formula.
\end{abstract}
\end{frontmatter}

 \section{Introduction }
\setcounter{equation}{0}
 In this paper, We consider  the following elliptic Hénon type equation
\begin{equation}\label{n}
\D^{2} u= |x|^a|u|^{p-1}u \; \,\,\mbox{in}\,\,\,\, \R^{n}_{+},  \quad \quad  u =\frac{\partial u}{\partial x_n} = 0 \quad  \mbox{on}\,\,\,\, \partial \R^{n}_{+},
\end{equation}
where $$p>1, \; \,\,n\geq 2,  \quad \; \,\,\R^n_+:=\{\; x\in \R^n,\;\; x_n>0\} \; \,\,\mbox{and}\,\,\,\, \partial \R^n_+:=\{ \; x\in \R^n,\;\; x_n=0\}.$$

\medskip

Our main objective is to classify the non-existence result for $C^4$-solutions for problems \eqref{n} belonging
to one of the following classes: stable solutions and solutions which are stable outside a compact set. 
\smallskip

We now list some known results. We start with the Liouville type theorems for the corresponding nonlinear problem

\begin{equation}\label{le}
 ( -\Delta)^m u  =|x|^{a}|u|^{p-1}u \,\, \mbox{in } \;\; \R^N,
 \end{equation}
have been largely studied in the literature ( see, e.g.,\cite{Farina, GidasSpruck, GS, DavilaDupaigneWangWei, hu, DancerDuGuo, WangYe}). In particular, when $m=1$ and $a=0,$  the first Liouville theorem was proved by Gidas and Spruck in \cite{GidasSpruck}, in which they proved that, for $1<p\leq\frac{n+2}{n-2}.$ Soon afterward, similar results were established in \cite{GS} for positive solutions of the
subcritical problem in the upper half-space 
\begin{equation*}
-\D u= |u|^{p-1}u \; \,\,\mbox{in}\,\,\,\, \R^{n}_{+},  \quad \quad  u =0 \quad  \mbox{on}\,\,\,\, \partial \R^{n}_{+},
\end{equation*}

After that, Chen and Li\cite{WX} obtained similar nonexistence results for the above two equations by using the moving plane method. In \cite{Farina}, Farina obtained the optimal Liouville type result for solutions stable at infinity. Indeed, he proved that a smooth nontrivial solution to \eqref{le} exists,
if and only if $p \geq p_{JL1}(n)$ and $n\geq 11,$ or $p=\frac{n+2}{n-2}$ and $n\geq 3.$ Here $p_{JL1}(n)$ denotes the so-called Joseph-Lundgren exponent (see \cite{CW, Farina}). In addition, similar results were established in \cite{Farina} for finite Morse index solutions in the upper half-space: $\R^{n}_{+},$ with homogeneous Dirichlet boundary conditions on $\partial \R^n_+$. Furthermore, strips provide an interesting example of unbounded domains where, as we shall see, rather sharp.

Furthermore, in a recent paper \cite{DancerDuGuo}, Dancer, Du and Guo extended some results in \cite{Farina} have considered \eqref{le} with $m=1$ and $a>-2,$. It was proved that there is no nontrivial stable solution  in $\R^n$ if $1<p<p_{JL1}(n,a)$ and that for $p\geq p_{JL1}(n,a),$ admits a positive radial stable solution in $\R^n$, where $p_{JL1}(n,a)$ is Joseph-Lundgren exponent for the Hénon type equation.
 In addition, Wang and Ye \cite{WangYe} obtained a Liouville-type result for finite Morse index solutions in $\R^n,$ which is a partial extension of results in \cite{DancerDuGuo}. 

 In a very interesting paper, D\'{a}vila et al \cite{DavilaDupaigneWangWei} investigated the bi-harmornic equation  i.e.~$m=2$ and $a=0,$ they derived a relevant monotonicity formula and employed blow down analysis to prove a sharp classification of stable at infinity solutions. However, for the fourth-order Hénon type equation  i.e.~$m=2$ and $a>0,$ studied by Hu \cite{hu}. He proved Liouville-type theorems for solutions belonging to one of the following classes: stable solutions and finite Morse index solutions (whether positive or sign-changing). His proof is based on a combination of the Pohozaev-type identity, monotonicity formula of solutions and a blowing down sequence. 
\medskip

Relying on  Hu’s approach \cite{hu} and using the  technics developed in \cite{Farina, DavilaDupaigneWangWei}, we give a Liouville-type theorems in the class of stable solution and finite Morse index solutions in the half space $\R^N_+$.  Before
stating our main results, we first recall the definition of such solutions.

\begin{defn}
We say that a solution $u$ of \eqref{n} belonging to $C^4(\overline{\R^n_+})$,

$\bullet$ is stable, if \begin{align}\label{quadr}
  Q_u (\psi):= \int_{\R^n_+}(\D \psi)^2 dx -p\int_{\R^n_+}|x|^a|u|^{P-1}\psi^2dx \geq 0,\quad
  \forall\,\psi \in C^2_c(\overline{\R^n_+}).
\end{align}

$\bullet$ is stable outside a compact set $\mathcal{K}\subset \R^n_+$, if $ Q_u (\psi)\geq0$ for any $\psi \in C^2_c(\overline{\R^n_+}\backslash \mathcal{K})$.

$\bullet$ More generally, the Morse index of a solution is defined as the maximal dimension of all subspaces $E$ of $C^2_c(\overline{\R^n_+})$ such that
$ Q_u (\psi)< 0$ in $E\backslash \{0\}$. Clearly, a solution stable if and only if its Morse index is equal to zero.
\end{defn}
\begin{re}
\begin{itemize}
\item[ (i).] Clearly a solution is stable if and only if its Morse index is equal to zero.
\item[ (ii).] Any finite Morse index solution $u$ is stable outside a compact set $\mathcal{K}\subset\R^n_+$. Indeed there exist $K\geq1$ and $X_K:=$ span$\{\varphi_1,\ldots,\varphi_K\}\subset C^2_c(\R^n_+)$ such that $Q_u(\varphi)<0$ for any $\varphi\in X_K\backslash\{0\}$. Then, $Q_u(\varphi)\geq0$ for every $\varphi\in C^2_c(\overline{\R^N_+}\backslash\mathcal{K})$, where $\mathcal{K}:=\cup^K_{j=1}\mbox{supp}(\varphi_j)$.
\end{itemize}
\end{re}

Now we can state our main results

\begin{thm}\label{th1.1}
Let $u\in C^{4}(\overline{\R^{n}_{+}})$ be a stable solution of \eqref{n}. If $1<p< p_{JL2}(n,a),$ then $u\equiv 0$
\end{thm}

\begin{thm}\label{th1.2}
Let $u \in C^{4}(\overline{\R^{n}_{+}})$ be a solution of \eqref{n} that is stable outside a compact set.\\
$\bullet$ If $1<p<p_{JL2}(n,0),\;\;p\neq \frac{n+4+2a}{n-4}$ then $u\equiv 0.$\\
$\bullet$ If $p=\frac{n+4+2a}{n-4},$ then $u$ has finite energy,i.e.,
$$\int_{\R^{n}_{+}}(\D u)^{2}=\int_{\R^{n}_{+}}|x|^a|u|^{p+1}<+\infty.$$
\end{thm}
\medskip
Here the representation of $p_{JL2}(n,a)$ in Theorem \ref{th1.1} is the fourth-order Joseph–Lundgren exponent which is computed by [6]. 

The organization of the rest of the paper is as follows. In section 1, we construct a monotonicity formula which is a crucial
tool to handle the supercritical case. In section 2, we establish some finer integral estimates
for the solutions of \eqref{n}. Also  we  obtain  a nonexistence result for the homogeneous stable solution of \eqref{n} in $\R_{+}^{n}\setminus \lbrace 0\rbrace,$ where $p$ belongs to $(\frac{n+4+2a}{n-4},p_{JL2}(n,a))$. Then we prove Liouville-type theorem for  stable solutions of $\eqref{n},$ this is Theorem \ref{th1.1} in section 3. To prove the result, we obtain some estimates of solutions, and show that the limit of blowing down sequence $u^{\infty}(x)=\displaystyle \lim_{\lambda\longrightarrow \infty}\lambda^{\frac{4+a}{p-1}}u(\lambda x)$ satisfies $E(u,r)\equiv const.$ Here, we use the monotonicity formula of ( Proposition 1.1 see below ). In section 4, we study Liouville-type theorem of finite Morse index solutions by the use of the Pohozaev-type identity, monotonicity formula and blowing down sequence. In the following, $C$
denotes always a generic positive constant, which could be changed from one line to another. 
\medskip

\subsection{\textbf{Monotonicity formula}}

\medskip
In this section, we construct a monotonicity formula which plays an important role in
dealing to understand supercritical elliptic equations or systems. This approach has been
used successfully for the Lane–Emden equation in \cite{DavilaDupaigneWangWei,hu}. Equation \eqref{n} has two important features. It is variational, with the energy functional given by
$$\int \bigg(\frac{1}{2}|\D u|^{2}-\frac{1}{p+1}|x|^a|u|^{p+1}\bigg).$$
For $\lambda>0,$  set $B^{+}_{\lambda}=B_{\lambda}\cap\R^{n}_{+}.$  Under the scaling transformation
$$u^{\lambda}(x)=\lambda^{\frac{4+a}{p-1}}u(\lambda x),$$
this suggests that the variations of the rescaled energy
$$\int_{B^{+}_{1}} \bigg(\frac{1}{2}|\D u^{\lambda}|^{2}-\frac{1}{p+1}|x|^a|u^{\lambda}|^{p+1}\bigg).$$
For any given $x \in \R^{n}_{+},$  we choose $u\in W^{4,2}_{loc}(\R^{n}_{+})\cap L^{p+1}_{loc}( \R^{n}_{+})$ and define
\begin{align*}
E(u,\lambda)&= \lambda^{ \frac{4(p+1)+2a}{p-1}-n}\left(\int_{B^{+}_{\lambda}}\frac{1}{2} (\D u)^{2}-\frac{1}{p+1}|x|^a|u|^{p+1} \right)
\\
&\;+  \frac{4+a}{2(p-1)}\left( n-2-\frac{4+a}{p-1}\right) \lambda^{\frac{8+2a}{p-1}+1-n}\int_{\partial B^{+}_{\lambda}}u^{2}\\
&\;+\frac{4+a}{2(p-1)}\left( n-2-\frac{4+a}{p-1}\right)\frac{d}{d\lambda}\left(\lambda^{\frac{8+2a}{p-1}+2-n}\int_{\partial B^{+}_{\lambda}}u^{2}\right)\\
&
+\frac{\lambda^{3}}{2}\frac{d}{d\lambda}\left[ \lambda^{\frac{8+2a}{p-1}+1-n}\int_{\partial B^{+}_{\lambda}}\left(\frac{4}{p-1}\lambda^{-1}u + \frac{\partial u}{\partial r}\right)^{2}\right]\\
&+\frac{1}{2}\frac{d}{d\lambda}\left[ \lambda^{\frac{8+2a}{p-1}+4-n}\int_{\partial B^{+}_{\lambda}}\left( |\nabla u |^{2}-\left|\frac{\partial u}{\partial r}\right|^{2} \right)\right]+ \frac{1}{2}\lambda^{\frac{8+2a}{p-1}+3-n}\int_{\partial B^{+}_{\lambda}}\left( |\nabla u |^{2}-\left|\frac{\partial u}{\partial r}\right|^{2} \right),
\end{align*}
where derivatives are taken in the sense of distributions.Then, we have the following monotonicity formula.
\begin{prop}\label{prop 2.1}
Assume that $n\geq5,\;a\geq0$ and $p>\frac{n+4+2a}{n-4},\;u\in W^{4,2}_{loc}(\R^{n}_{+})$ and $|x|^a|u|^{p+1} \in L^{1}_{loc}( \R^{n}_{+})$ be a weak solution of \eqref{n}. Then, $E(u,\lambda)$ is non-decreasing in $\lambda>0.$ Furthermore there is a constant $C(n,p,a)>0$ such that
\begin{align}
\frac{d}{dr}E(u,\lambda)\geq C(n,p,a)\lambda^{-n+2+\frac{8+2a}{p-1}}\int_{\partial B^{+}_{\lambda}}\left(\frac{4+a}{p-1}\lambda^{-1}u + \frac{\partial u}{\partial r}\right)^{2}dS.
\end{align}
\end{prop}

{\bf Proof.} The proof follows the main lines of the demonstration of Theorem  2.1 in \cite{hu}, with
small modifications. Since the boundary integrals in $E(u,\lambda)$ only involve second order derivatives of $u,$ the boundary integrals in $\frac{dE}{d\lambda}(u,\lambda)$ only involve third order derivatives of $u.$  Thus, the following calculations can be rigorously verified. Assume that $x=0$ and that the balls $B_{\lambda}$ are all centered at $0.$ Take

 $$\tilde{E}(\lambda)= \lambda^{ \frac{4(p+1)+2a}{p-1}-n}\int_{B^{+}_{\lambda}}\frac{1}{2} (\D u)^{2}-\frac{1}{p+1}|x|^a|u|^{p+1} .$$
 Define $$v=\D u, \quad u^{\lambda}(x)=\lambda^{\frac{4+a}{p-1}}u(\lambda x) \quad\mbox{and}\quad v^{\lambda}(x)=\lambda^{\frac{4+a}{p-1}+2}v(\lambda x).$$
We still have $v^{\lambda}=\D u^{\lambda},\; \D v^{\lambda}=|x|^a|u^{\lambda}|^{p-1}u^{\lambda},$ and by differentiating in $\lambda,$
$$\D \frac{du^{\lambda}}{d\lambda}=\frac{dv^{\lambda}}{d\lambda}.$$
Note that differentiation in $\lambda$ commutes with differentiation and integration in $x.$ A rescaling shows
 $$\tilde{E}(\lambda)= \int_{B^{+}_{1}}\frac{1}{2} (v^{\lambda})^{2}-\frac{1}{p+1}|x|^a|u^{\lambda}|^{p+1} .$$
hence
\begin{align*}
\frac{d}{d\lambda}\tilde{E}(\lambda)&=\int_{B^{+}_{1}}v^{\lambda}\frac{dv^{\lambda}}{d\lambda}-|x|^a|u^{\lambda}|^{p-1}u^{\lambda}\frac{du^{\lambda}}{d\lambda}\\
&\;=\int_{B^{+}_{1}}v^{\lambda}\D \frac{du^{\lambda}}{d\lambda}-\D v^{\lambda}\frac{du^{\lambda}}{d\lambda}=\int_{\partial B^{+}_{1}}v^{\lambda}\frac{\partial}{\partial r}\frac{du^{\lambda}}{d\lambda}-\frac{\partial v^{\lambda}}{\partial r}\frac{du^{\lambda}}{d\lambda}.
\end{align*}

Since $u^{\lambda} = 0 $ in $\partial \R^n_+$ for any $\lambda >0,$ then  $\frac{d u^{\lambda}}{d\lambda}=0$ in $\partial \R^n_+.$ Hence, all boundary terms appearing in the integrations by parts vanish under the  Dirichlet boundary conditions. So, we get
\begin{align}\label{2.4}
\frac{d}{d\lambda}\tilde{E}(\lambda)=\int_{\partial B^{+}_{1}}\left(v^{\lambda}\frac{\partial}{\partial r}\frac{du^{\lambda}}{d\lambda}-\frac{\partial v^{\lambda}}{\partial r}\frac{du^{\lambda}}{d\lambda}\right).
\end{align}
In what follows, we express all derivatives of $u^{\lambda}$ in the $r=|x|$ variable in terms of derivatives in the $\lambda$ variable. In the definition of $u^{\lambda}$ and $v^{\lambda}$, directly differentiating in $\lambda$ gives
\begin{align}\label{2.5}
\frac{du^{\lambda}}{d\lambda}(x)=\frac{1}{\lambda}\left(\frac{4+a}{p-1}u^{\lambda}(x)+r\frac{\partial u^{\lambda}}{\partial r}(x)\right),
\end{align}
and
\begin{align}\label{2.6}
\frac{dv^{\lambda}}{d\lambda}(x)=\frac{1}{\lambda}\left(\frac{2(p+1)+a}{p-1}v^{\lambda}(x)+r\frac{\partial v^{\lambda}}{\partial r}(x)\right).
\end{align}
In \eqref{2.5}, taking derivatives in $\lambda$ once again, we get
\begin{align}\label{2.7}
\lambda\frac{d^{2}u^{\lambda}}{d\lambda^{2}}(x)+\frac{du^{\lambda}}{d\lambda}(x)=\frac{4+a}{p-1}\frac{du^{\lambda}}{d\lambda}(x)+r\frac{\partial}{\partial r}\frac{du^{\lambda}}{d\lambda}(x).
\end{align}
Substituting \eqref{2.6} and \eqref{2.7} into \eqref{2.4}, we obtain
\begin{align}
\begin{split}
\label{2.8}
\frac{d\tilde{E}}{d\lambda}&=\int_{\partial B^{+}_{1}}v^{\lambda}\left(\lambda\frac{d^{2}u^{\lambda}}{d\lambda^{2}}+\frac{p-5-a}{p-1}\frac{du^{\lambda}}{d\lambda}\right)-\frac{du^{\lambda}}{d\lambda}\left(\lambda\frac{dv^{\lambda}}{d\lambda}-\frac{2(p+1)+a}{p-1}v^{\lambda}\right)
\\&=\int_{\partial B^{+}_{1}}\lambda v^{\lambda}\frac{d^{2}u^{\lambda}}{d\lambda^{2}} + 3v^{\lambda}\frac{du^{\lambda}}{d\lambda}-\lambda\frac{du^{\lambda}}{d\lambda}\frac{dv^{\lambda}}{d\lambda}.
\end{split}
\end{align}
Observe that $v^{\lambda}$ is expressed as a combination of $x$ derivatives of $u^{\lambda}$. So we also transform $v^{\lambda}$ into $\lambda$ derivatives of $u^{\lambda}$. By taking derivatives in $r $ in \eqref{2.5} and noting \eqref{2.7}, we get on $\partial B^{+}_{1},$
\begin{eqnarray*}
\frac{\partial^{2}u^{\lambda}}{\partial r^{2}}&=& \lambda \frac{\partial}{\partial r}\frac{\partial u^{\lambda}}{\partial \lambda}-\frac{p+3+a}{p-1}\frac{\partial u^{\lambda}}{\partial r}\nonumber\\&=&\lambda^{2} \frac{\partial^{2} u^{\lambda}}{\partial \lambda^{2}}+\frac{p-5-a}{p-1} \lambda\frac{d u^{\lambda}}{d \lambda}-\frac{p+3+a}{p-1}\left(\lambda\frac{d u^{\lambda}}{d \lambda}-\frac{4+a}{p-1}u^{\lambda}\right)\nonumber\\&=&\lambda^{2} \frac{\partial^{2} u^{\lambda}}{\partial \lambda^{2}}-\frac{8+2a}{p-1}\lambda\frac{d u^{\lambda}}{d \lambda}+\frac{(4+a)(p+3+a)}{(p-1)^{2}}u^{\lambda}.
\end{eqnarray*}
Then on $\partial B^{+}_{1},$
\begin{eqnarray*}
v^{\lambda}&=&\frac{\partial^{2}u^{\lambda}}{\partial r^{2}}+\frac{n-1}{r}\frac{\partial u^{\lambda}}{\partial r}+\frac{1}{r^{2}}\D_{\theta} u^{\lambda}\nonumber\\&=&\lambda^{2} \frac{d^{2} u^{\lambda}}{d \lambda^{2}}-\frac{8+2a}{p-1}\lambda\frac{d u^{\lambda}}{d \lambda}+\frac{(4+a)(p+3+a)}{(p-1)^{2}}u^{\lambda}+(n-1)\left(\lambda\frac{d u^{\lambda}}{d \lambda}-\frac{4+a}{p-1}u^{\lambda}\right)+\D_{\theta} u^{\lambda}\\&=&\lambda^{2} \frac{d^{2} u^{\lambda}}{d \lambda^{2}}+\left(n-1-\frac{8+2a}{p-1}\right)\lambda\frac{d u^{\lambda}}{d \lambda}+\frac{4+a}{p-1}\left(\frac{4+a}{p-1}-n+2\right)u^{\lambda}+\D_{\theta} u^{\lambda}.
\end{eqnarray*}
Here $\D_{\theta}$ is the Beltrami–Laplace operator on $\partial B_{1}$ and below $\nabla_{\theta}$ represents the tangential derivative on  $\partial B_{1}$. For notational convenience, we also define the constants

$$\alpha=n-1-\frac{8+2a}{p-1},\;\;\;\beta=\frac{4+a}{p-1}\left(\frac{4+a}{p-1}-n+2\right).$$
Now \eqref{2.8} reads
\begin{align*}
\frac{d}{d\lambda}\tilde{E}(\lambda):=I_{1}+I_{2},
\end{align*}
where
\begin{align*}
I_{1}&:= \int_{\partial B^{+}_{1}}\lambda\left(\lambda^{2}\frac{d^{2} u^{\lambda}}{d \lambda^{2}}+\alpha\lambda\frac{d u^{\lambda}}{d \lambda}+\beta u^{\lambda}\right)\frac{d^{2} u^{\lambda}}{d \lambda^{2}}\\
&+3\left(\lambda^{2}\frac{d^{2} u^{\lambda}}{d \lambda^{2}}+\alpha\lambda\frac{d u^{\lambda}}{d \lambda}+\beta u^{\lambda}\right)\frac{d u^{\lambda}}{d \lambda}-\lambda\frac{d u^{\lambda}}{d \lambda}\frac{d}{d \lambda}\left(\lambda^{2}\frac{d^{2} u^{\lambda}}{d \lambda^{2}}+\alpha\lambda\frac{d u^{\lambda}}{d \lambda}+\beta u^{\lambda}\right),
\end{align*}
and $$I_{2}:=\int_{\partial B^{+}_{1}}\lambda\D_{\theta}u^{\lambda}\frac{d^{2} u^{\lambda}}{d \lambda^{2}}+3\D_{\theta}u^{\lambda}\frac{d u^{\lambda}}{d \lambda}-\lambda\frac{d u^{\lambda}}{d \lambda}\D_{\theta}\frac{d u^{\lambda}}{d \lambda}.$$

 Let  $\lambda >0.$ Since  $\frac{d u^{\lambda}}{d\lambda}=0$ in $\partial \R^n_+$ then, all boundary terms appearing in the integrations by parts vanish under the  Dirichlet boundary conditions, hence the calculations are even
easier. The integral  $I_{2}$ can be estimated as

\begin{eqnarray*}
I_{2}&=&\int_{\partial B^{+}_{1}}-\lambda\nabla_{\theta}u^{\lambda}\nabla_{\theta}\frac{d^{2} u^{\lambda}}{d \lambda^{2}}-3\nabla_{\theta}u^{\lambda}\nabla_{\theta}\frac{d u^{\lambda}}{d \lambda}+\lambda\left|\nabla_{\theta}\frac{d u^{\lambda}}{d \lambda}\right|^{2}\nonumber\\&=&-\frac{\lambda}{2}\frac{d^{2}}{d\lambda^{2}}\left(\int_{\partial B^{+}_{1}}|\nabla_{\theta} u^{\lambda}|^{2}\right)-\frac{3}{2}\frac{d}{d\lambda}\left(\int_{\partial B^{+}_{1}}|\nabla_{\theta} u^{\lambda}|^{2}\right)+2\lambda\int_{\partial B^{+}_{1}}\left|\nabla_{\theta}\frac{d u^{\lambda}}{d \lambda}\right|^{2}\nonumber\\&=&-\frac{1}{2}\frac{d^{2}}{d\lambda^{2}}\left(\lambda \int_{\partial B^{+}_{1}}|\nabla_{\theta} u^{\lambda}|^{2}\right)-\frac{1}{2}\frac{d}{d\lambda}\left(\int_{\partial B^{+}_{1}}|\nabla_{\theta} u^{\lambda}|^{2}\right)+2\lambda\int_{\partial B^{+}_{1}}\left|\nabla_{\theta}\frac{d u^{\lambda}}{d \lambda}\right|^{2}\nonumber\\&\geq&-\frac{1}{2}\frac{d^{2}}{d\lambda^{2}}\left(\lambda \int_{\partial B^{+}_{1}}|\nabla_{\theta} u^{\lambda}|^{2}\right)-\frac{1}{2}\frac{d}{d\lambda}\left(\int_{\partial B^{+}_{1}}|\nabla_{\theta} u^{\lambda}|^{2}\right).
\end{eqnarray*}
Furthermore, a direct calculation implies that
\begin{align*}
I_{1}&=  \int_{\partial B^{+}_{1}}\lambda^{3}\left(\frac{d^{2} u^{\lambda}}{d \lambda^{2}}\right)^{2}+\lambda^{2}\frac{d^{2} u^{\lambda}}{d \lambda^{2}}\frac{d u^{\lambda}}{d \lambda}+\beta \lambda u^{\lambda}\frac{d^{2} u^{\lambda}}{d \lambda^{2}}+3\beta u^{\lambda}\frac{d u^{\lambda}}{d \lambda}+(2\alpha-\beta)\lambda\left(\frac{d u^{\lambda}}{d \lambda}\right)^{2}-\lambda^{3}\frac{d u^{\lambda}}{d \lambda}\frac{d^{3} u^{\lambda}}{d \lambda^{3}}\\
&= \int_{\partial B^{+}_{1}}2\lambda^{3}\left(\frac{d^{2} u^{\lambda}}{d \lambda^{2}}\right)^{2}+4\lambda^{2}\frac{d^{2} u^{\lambda}}{d \lambda^{2}}\frac{d u^{\lambda}}{d \lambda}
+(2\alpha-2\beta)\lambda\left(\frac{d u^{\lambda}}{d \lambda}\right)^{2}+\frac{\beta}{2}\frac{d^{2}}{d\lambda^{2}}[\lambda(u^{\lambda})^{2}]+\frac{\beta}{2}\frac{d}{d\lambda}(u^{\lambda})^{2}\\
&-\frac{1}{2}\frac{d}{d\lambda}\left[\lambda^{3}\frac{d}{d\lambda}\left(\frac{d u^{\lambda}}{d \lambda}\right)^{2}\right].
\end{align*}
Here we have used the relations (writing $f^{'}=\frac{d}{d\lambda}f$ etc.)
$$\lambda ff^{''}=\left(\frac{\lambda}{2}f^{2}\right)^{''}-2ff^{'}-\lambda(f^{'})^{2},\quad\mbox{and}\quad-\lambda^{3}f^{'}f^{'''}=-\left[\frac{\lambda^{3}}{2}((f^{'})^{2})^{'}\right]^{'}+3\lambda^{2}f^{'}f^{''}+\lambda^{3}(f^{''})^{2}.$$
Since $p>\frac{n+4+2a}{n-4},$ direct calculations show that
\begin{eqnarray}
\alpha-\beta=\left(n-1-\frac{8+2a}{p-1}\right)-\frac{4+a}{p-1}\left(\frac{4+a}{p-1}-n+2\right)>1.
\end{eqnarray}
Consequently,
\begin{align*}
&\;2\lambda^{3}\left(\frac{d^{2} u^{\lambda}}{d \lambda^{2}}\right)^{2}+4\lambda^{2}\frac{d^{2} u^{\lambda}}{d \lambda^{2}}\frac{d u^{\lambda}}{d \lambda}
+(2\alpha-2\beta)\lambda\left(\frac{d u^{\lambda}}{d \lambda}\right)^{2}
\\
&\;=2\lambda\left(\lambda\frac{d^{2} u^{\lambda}}{d \lambda^{2}}+\frac{d u^{\lambda}}{d \lambda}\right)^{2}+(2\alpha-2\beta-2)\lambda\left(\frac{d u^{\lambda}}{d \lambda}\right)^{2}\geq0.
\end{align*}
We we conclude then $$I_{1}\geq \int_{\partial B^{+}_{1}}\frac{\beta}{2}\frac{d^{2}}{d\lambda^{2}}[\lambda(u^{\lambda})^{2}]-\frac{1}{2}\frac{d}{d\lambda}\left[\lambda^{3}\frac{d}{d\lambda}\left(\frac{d u^{\lambda}}{d \lambda}\right)^{2}\right]+\frac{\beta}{2}\frac{d}{d\lambda}(u^{\lambda})^{2}.$$
Now, rescaling back, we can write those $\lambda$ derivatives in $I_{1}$ and $I_{2}$ as follows.
$$\int_{\partial B^{+}_{1}}\frac{d}{d\lambda}(u^{\lambda})^{2}=\frac{d}{d\lambda}\left(\lambda^{\frac{8+2a}{p-1}+1-n}\int_{\partial B^{+}_{\lambda}} u^{2}\right),$$
$$\int_{\partial B^{+}_{1}} \frac{d^{2}}{d \lambda^{2}}[\lambda(u^{\lambda})^{2}]=\frac{d^{2}}{d\lambda^{2}}\left(\lambda^{\frac{8+2a}{p-1}+2-n}\int_{\partial B^{+}_{\lambda}} u^{2}\right),$$
$$\int_{\partial B^{+}_{1}}\frac{d}{d\lambda}\left[\lambda^{3}\frac{d}{d\lambda}\left(\frac{d u^{\lambda}}{d \lambda}\right)^{2}\right]=\frac{d}{d\lambda}\left[\lambda^{3}\frac{d}{d\lambda}\left(\lambda^{\frac{8+2a}{p-1}+1-n}\int_{\partial B^{+}_{\lambda}} \left( \frac{4+a}{p-1} \lambda^{-1} u +\frac{\partial u}{\partial r} \right)^{2} \right) \right] ,$$
$$\frac{d^{2}}{d\lambda^{2}} \left( \lambda \int_{\partial B^{+}_{1}} |\nabla_{\theta} u^{\lambda}|^{2} \right)= \frac{d^{2}}{d\lambda^{2}}\left[\lambda^{1+\frac{8+2a}{p-1}+2+1-n}\int_{\partial B^{+}_{\lambda}}\left( |\nabla u |^{2}-\left|\frac{\partial u}{\partial r}\right|^{2} \right)\right] , $$
and
$$\frac{d}{d\lambda}\left(\int_{\partial B^{+}_{1}}|\nabla_{\theta}u^{\lambda}|^{2} \right)=\frac{d}{d\lambda}\left[\lambda^{\frac{8+2a}{p-1}+2+1-n}\int_{\partial B^{+}_{\lambda}}\left( |\nabla u |^{2}-\left|\frac{\partial u}{\partial r}\right|^{2} \right)\right].$$

Substituting these into $\frac{d}{d\lambda} E(u,\lambda)$ we finish the proof.\qed

\medskip
\section{Main technical tool.}
\setcounter{equation}{0}
\medskip

\subsection{\textbf{Integral estimates}}

\medskip

For $\beta>0,$  set  $B^{+}_{\beta}=B_{\beta}\cap\R^{n}_{+}$ and $ A^{+}_{\beta}=\{x \in \R^{n}_{+},  a_{1}\beta <|x|<a_{2}\beta\},$  for some $0<a_{1}<a_{2}.$ Let $u$ be a solution of \eqref{n}, which is stable outside a compact set $\mathcal{K} \subset B^{+}_{R_0}.$  For all $R> 4R_0 $, we define a family of test functions $  \psi= \psi _{(R,R_0) }\in C^{2}_c(\R^N)$ satisfying

\begin{align}\label{eq1}
 \left\{\begin{array}{llllllllllllllllll}
 0\leq \psi \leq 1 \mbox{ and } \psi\equiv 0\;\quad \mbox{if }\quad |x|<R_0  \;\mbox{or }\; |x|> 2 R,\\
\psi\equiv 1 \;\quad \mbox{if } \quad 2R_0<|x|< R ,\\\\
|\nabla^q \psi|\leq C R_0^{-q} \;  \quad \mbox{if }\;\;\; R_0 <|x|<2R_0, \\
 |\nabla^q \psi|\leq CR^{-q} \quad \; \mbox{if } \quad  R<|x|<2R,\; \mbox{and } 1\leq q\leq 4.
\end{array}
\right.
\end{align}

Similarly, if $u$ is a stable  solution of \eqref{n}, then $\psi=\psi_{(R)}$, with $R>0$  verifying \eqref{eq1} with $R_0=0$ that is $\psi=1 \mbox{ if } |x|< R.$
\medskip

First of all, we need the following lemma which plays an important role in dealing with  Theorem \ref{th1.1} and  Theorem \ref{th1.2}

\begin{lem} \label{2.3}
 Let $u\in C^4(\overline{\R^n_+})$ be a solution of \eqref{n}, which is stable outside a compact set $\mathcal{K}.$ Let $R_0>0$ such that $\mathcal{K} \subset B^{+}_{R_0}$ and set  $v=\D u$ , there hold
\begin{align}\label{2.12}
\int_{B^{+}_{R}}v^2  +\int_{B^{+}_{R}}|x|^a|u|^{p+1}\leq C_0+ C R^{-4}\int_{A^{+}_{R}}u^{2}+ C R^{-2}\int_{A^{+}_{R}}|uv|,\;\quad \forall\;\; R>4R_0,
\end{align}
and
\begin{align}\label{2.13}
\int_{B^{+}_{R}}v^{2}+ \int_{B^{+}_{R}}|x|^a|u|^{p+1} \leq C(1+ R^{n-\frac{4(p+1)+2a}{p-1}}),\;\quad \forall\;\; R>4R_0.
\end{align}
\end{lem}
\noindent{\bf Proof.}
\smallskip

 \noindent{\bf Proof of \eqref{2.13}.} First, for $\epsilon \in (0, 1)$ and $\eta \in C^2(\R^N)$,
\begin{align*}
\int_{\R^{n}_{+}} [\Delta (u\eta)]^2  & = \int_{\R^{n}_{+}}\left(u \D \eta+2\nabla u\nabla \eta + \eta\D u\right)^{2} \\
& \leq \left(1+C\epsilon\right)\int_{\R^{n}_{+}}v^2 {\eta}^2 + \frac{C}{\e}\int_{\R^{n}_{+}}u^2
  (\D\eta)^2  +\frac{C}{\epsilon} \int_{\R^{n}_{+}} |\nabla u|^2|\nabla \eta|^2 .
   \end{align*}
  Using $\Delta(u^2)  = 2|\nabla u|^2  + 2u\D u$,
\begin{align}
\label{newest5}
2\int_{\R^{n}_{+}}|\nabla u|^{2}|\nabla\eta|^{2}= \int_{\R^{n}_{+}} u^{2}\D(|\nabla\eta|^{2})-2\int_{\R^{n}_{+}}uv|\nabla\eta|^{2}.
\end{align}

So, we we get
\begin{align}
\label{newestj5}
\int_{\R^{n}_{+}}[\Delta (u\eta)]^2\leq  \left(1+C\epsilon\right)\int_{\R^{n}_{+}}v^2 \eta^{2}+ C_{\epsilon}\int_{\R^{n}_{+}}u^{2}\Big[(\D\eta)^{2}+|\D(|\nabla\eta|^{2})|\Big]+ \frac{C}{\e} \int_{\R^{n}_{+}}|uv||\nabla\eta|^{2} .
\end{align}
 Take $\eta = \eta^m$ with $m \geq 2$. Apply Cauchy-Schwarz's inequality, we get
  \begin{align}
  \label{t43}
  \int_{\R^{n}_{+}}|uv||\nabla \eta^{m}|^2 \leq C\epsilon^2 \int_{\R^{n}_{+}}v^2 \eta^{2m} +
  C_{\epsilon, m}\int_{\R^{n}_{+}} u^2 |\nabla \eta|^4\eta^{2m-4}.
\end{align}
Substitute $\eta$ by $\psi^{m}$ in \eqref{newestj5}, then from \eqref{t43} and \eqref{eq1}, we obtain

$$\int_{B^{+}_{2R}}[\Delta (u\psi^{m})]^2\leq C_0+ \left(1+C\epsilon\right)\int_{B^{+}_{2R}}v^2 \psi^{2m}+ C_{\epsilon} R^{-4}\int_{A^{+}_{R}}u^{2},$$
where $$C_0=C R^{-4}\int_{A^{+}_{0}}u^{2} ,\quad \; \mbox{and } \quad   A^{+}_{0}=\{x \in \R^{n}_{+},\;\;  R_0 <|x|<2R_0\}.$$

Let $u$ be a solution of \eqref{n}, which is stable outside a compact set $\mathcal{K} \subset B^{+}_{R_0}.$  Clearly $u\psi^{m}\in H^{2}_{0}(B^{+}_{2R}\setminus B^{+}_{R_{0}})$ so after a standard approximation argument, the main inequality of stability \eqref{quadr} implies

$$p\int_{B^{+}_{2R}}|x|^a|u|^{p+1}\psi^{2m} -\int_{B^{+}_{2R}}(\D(u\psi^{m}))^{2}\leq0,\;\quad \forall\;\; R>4R_0.$$
Therefore, we conclude then
\begin{align}
\label{newest6}
p\int_{B^{+}_{2R}}|x|^a|u|^{p+1}\psi^{2m}-\left(1+C\epsilon\right)\int_{B^{+}_{2R}}v^2 \psi^{2m}\leq C_0+ C_{\epsilon} R^{-4}\int_{A^{+}_{R}}u^{2}.
\end{align}

On the other hand, recall that $u =\frac{\partial u}{\partial x_n} = 0 \,\  \mbox{in}\,\, \partial \R^{n}_{+}.$ Multiply the equation \eqref{n} by $u\eta^2$, $\eta\in C^2(\R^N)$ and
integrate by parts, using again \eqref{newest5}
\begin{align}
\label{ewe}
  \begin{split}
& \; \int_{\R^{n}_{+}}\Big[v^2 \eta ^2 - |x|^a|u|^{p+1} \eta^2\Big] \\
= & \; -4\int_{\R^{n}_{+}}\eta v \nabla u\cdot\nabla \eta -2\int_{\R^{n}_{+}}\eta u v \D \eta  - 2\int_{\R^{n}_{+}}uv|\nabla \eta|^2\\
\leq & \; C\e\int_{\R^{n}_{+}} v^2 \eta ^2  + C_\e\int_{\R^{n}_{+}} u^2(\Delta \eta)^2  + C_\e \int_{\R^{n}_{+}} |\nabla u|^2|\nabla\eta|^2 - 2\int_{\R^{n}_{+}}uv|\nabla \eta|^2 \\
\leq & \; C\e\int_{\R^{n}_{+}} v^2 \eta ^2  + C_\e\int_{\R^{n}_{+}} u^2\Big[(\Delta \eta)^2 + |\Delta(|\nabla \eta|^2)\Big]  + C_\e\int_{\R^{n}_{+}}|uv||\nabla \eta|^2.
 \end{split}
\end{align}

Using the above inequality (where one substitutes $\eta$ by $\psi^{m}$),  it follows from \eqref{t43} and  \eqref{eq1} that
\begin{align}
\label{newest7}
\left(1-C\epsilon\right)\int_{B^{+}_{2R}}v^2 \psi^{2m}-\int_{B^{+}_{2R}}|x|^a|u|^{p+1}\psi^{2m}\leq C_0+ C_{\epsilon} R^{-4}\int_{A^{+}_{R}}u^{2}.
\end{align}

Taking $\epsilon > 0$ but small enough, multiplying \eqref{newest7} by $\frac{1+2C\epsilon}{1-C\epsilon}$, adding it with \eqref{newest6}  we get then
  \begin{align*}
C\epsilon\int_{B^{+}_{2R}}v^2 \psi^{2m} + \left(p-\frac{1+2C\epsilon}{1-C\epsilon}\right)\int_{B^{+}_{2R}}|x|^a|u|^{p+1}\psi^{2m}\leq C_0+ C_{\epsilon} R^{-4}\int_{A^{+}_{R}}u^{2}.
\end{align*}
As $p > 1$ and $A^{+}_{R}\subset B^{+}_{2R} ,$  using $\epsilon > 0$ small enough, there holds
\begin{align*}
\int_{B^{+}_{R}}v^2  +\int_{B^{+}_{R}}|x|^a|u|^{p+1}\leq C_0+ C R^{-4}\int_{A^{+}_{R}}u^{2}.
\end{align*}
Applying Young's inequality, we deduce then  for any $\epsilon'>0$
\begin{align*}
\int_{B^{+}_{R}}v^2  +\left(1-\epsilon'\right)\int_{B^{+}_{R}}|x|^a|u|^{p+1}\leq C_0+ CR^{n-\frac{4(p+1)+2a}{p-1}},\;\quad \forall\;\; R>4R_0.
\end{align*}
Take $\epsilon' > 0$ small enough, the estimate \eqref{2.13} is proved.

\smallskip

Let-us now prove \eqref{2.12}.  Invoking now \eqref{newestj5} where we substitute $\eta $ by $\psi^m,$ we obtain
 $$\int_{B^{+}_{2R}}[\Delta (u\psi^{m})]^2\leq C_0+ \left(1+C\epsilon\right)\int_{B^{+}_{2R}}v^2 \psi^{2m}+ C_{\epsilon} R^{-4}\int_{A^{+}_{R}}u^{2} +C_{\epsilon} R^{-2}\int_{A^{+}_{R}}|uv|.$$

Adopting the similar argument as  above where we use  the equality \eqref{ewe} and inequality of stability \eqref{quadr}, we obtain  readily the estimate \eqref{2.12}. Thus, Lemma \ref{2.3} is well proved.\qed

\medskip
\subsection {\textbf{Homogeneous solutions}}
\medskip

In this section, we obtain a nonexistence result for a homogeneous stable solution of \eqref{n}.
\begin{prop}\label{prop1}
Let $u\in W^{2,2}_{loc}(\R^n_{+}\backslash \{0\})$ be a homogeneous, stable solution of \eqref{n} in $\R^n_{+}\backslash \{0\}$, $p \in  \left(\frac{n+4+2a}{n-4},\; p_{JL2}(n,a)\right)$. Assume that $|x|^a|u|^{p+1}\in L^{1}_{loc}(\R^n_{+} \backslash \{0\})$. Then $u\equiv 0$.
\end{prop}
\noindent{\bf Proof.} Let $u$ be a homogeneous solution of \eqref{n}, that is there exists a $w \in W^{2,2}(\mathbb{S}^{n-1}_{+})$ such that in polar coordinates $$u(r, \theta) =  r^{-\frac{4+a}{p-1}} w (\theta).$$
 Denote $A^{+}_R= B^{+}_{2R}\backslash B^{+}_R.$ Since $u \in W^{2,2}(A^{+}_1)$ and $|x|^a|u|^{p+1}\in L^{1}(A^{+}_1)$, it implies that $$w \in W^{2,2}(\mathbb{S}^{n-1}_{+}) \cap L^{p+1}(\mathbb{S}^{n-1}_{+}).$$
A direct calculation gives
\begin{align} \label{e:022}
\D^2_{\theta} w(\theta)- J_1 \D_{\theta} w(\theta) + J_2 w(\theta)= |w|^{p-1} w\;\,\,\,\, \,\,\mbox{in}\,\,\,\,\; \mathbb{S}^{n-1}_{+}, \quad \quad w =\frac{\partial w}{\partial \theta_n}= 0 \,\,\,\,\mbox{on}\,\,\,\, \partial\mathbb{S}^{n-1}_{+},
\end{align}
where
 $$J_1= \left(\frac{4+a}{p-1} +2\right)\left(n-4-\frac{4+a}{p-1}\right)+  \frac{4+a}{p-1} \left(n-2-\frac{4+a}{p-1}\right),$$
 and
$$J_2=\frac{4+a}{p-1} \left(\frac{4+a}{p-1} +2\right)\left(n-4-\frac{4+a}{p-1} \right)\left(n-2-\frac{4+a}{p-1} \right). $$
 Because $w \in W^{2,2}(\mathbb{S}^{n-1}_{+})$, we can test \eqref{e:022} with $w$, and we obtain
\begin{eqnarray} \label{e:023}
\int_{\mathbb{S}^{n-1}_{+}} (\D _{\theta} w)^2 + J_1 |\nabla _{\theta}w|^2 + J_2 w^2\;d\theta = \int_{\mathbb{S}^{n-1}_{+}} |w|^{p+1}d\theta.
 \end{eqnarray}
As in \cite{DavilaDupaigneWangWei}, for any $\epsilon > 0 $, choose an $\eta_{\epsilon}\in C^{\infty}_0 ((\frac{\epsilon}{2}, \, \frac{2}{\epsilon})) $, such that $\eta_{\epsilon}\equiv 1$ in $(\epsilon, \, \frac{1}{\epsilon})$, and $$ r|\eta'_{\epsilon}(r)|+ r^2 |\eta''_{\epsilon}(r)|\leq 64, \quad for\;all \; r>0. $$
We assume that $\Omega_{k}=B_{2k/\epsilon}\backslash B_{\epsilon/2k}.$  Since $ w \in W ^{2,2}(\mathbb{S}^{n-1}_{+})\cap  L^{p+1}(\mathbb{S}^{n-1}_{+})$,\,\, $r^{-\frac{n-4}{2}} w(\theta) \eta_{\epsilon}(r)$ can be approximated by  $C^{\infty}_0(\Omega_{2}\cap\R^n_{+})$ functions in $W^{2, \,2}(\Omega_{1}\cap\R^n_{+} )\cap L^{p+1} (\Omega_{1}\cap\R^n_{+} )$. Hence, in the stability condition for $u$, we are allowed to choose a test function of the form $$r^{-\frac{n-4}{2}} w(\theta)\eta_{\epsilon}(r).$$
Direct calculations show that
\begin{eqnarray*}
\D \left(r^{-\frac{n-4}{2}} w(\theta)\eta_{\epsilon}(r)  \right)&=& -\frac{n(n-4)}{4} r^{-\frac{n}{2}} w(\theta) \eta_{\epsilon}(r)+ 3 r^{-\frac{n}{2}+1} w(\theta) \eta'_{\epsilon}(r)\nonumber\\&+& r^{-\frac{n}{2}+2}w(\theta) \eta''_{\epsilon}(r) + r^{-\frac{n}{2}} \D_{\theta} w(\theta) \eta_{\epsilon}(r),
\end{eqnarray*}

Substituting this into the stability condition for $u$, we deduce that
\begin{align*}
&p \left(\int_{\mathbb{S}^{n-1}_{+}} |w|^{p+1} d\theta \right)\left( \int_0^{+\infty} r^{-1} \eta_{\epsilon}(r)^2 dr \right)
 \\
\leq & \;\left(\int_{\mathbb{S}^{n-1}_{+}}\left((\D _{\theta} w)^2+\frac{n(n-4)}{2}|\nabla _{\theta}w|^2+\frac{n^2(n-4)^2}{16}w^2\right) d\theta\right)\left( \int_0^{+\infty} r^{-1} \eta_{\epsilon}(r)^2 dr \right)
\\
+ & \;O \bigg [\int_0^{+\infty} \left( r \eta'_{\epsilon}(r)^2 + r^3 \eta''_{\epsilon}(r)^2 +\eta_{\epsilon}(r)|\eta'_{\epsilon}(r)|+r \eta_{\epsilon}(r)|\eta''_{\epsilon}(r)| \right) dr \times \int_{\mathbb{S}^{n-1}_{+}}\left( |\nabla _{\theta}w(\theta)|^2 + w(\theta)^2\right)d\theta \bigg].
\end{align*}
Note that $$\int_0^{+\infty} r^{-1} \eta_{\epsilon}(r)^2 dr\geq |\log \epsilon |,$$
$$\int_0^{+\infty} \left( r \eta'_{\epsilon}(r)^2 + r^3 \eta''_{\epsilon}(r)^2 +\eta_{\epsilon}(r)|\eta'_{\epsilon}(r)| +r \eta_{\epsilon}(r)|\eta''_{\epsilon}(r)| \right) dr\leq C,$$
for some constant $C$ independent of $\epsilon$. By letting $\epsilon\rightarrow 0$, we obtain
\begin{eqnarray*}
p \int_{\mathbb{S}^{n-1}_{+}} |w|^{p+1} d\theta \leq \int_{\mathbb{S}^{n-1}_{+}} (\D _{\theta} w)^2 +\frac{n(n-4)}{2} |\nabla _{\theta}w|^2 +\frac{n^2(n-4)^2}{16}w^2 \; d\theta.
\end{eqnarray*}
Substituting \eqref{e:023} into this we derive
\begin{align*}
\int_{\mathbb{S}^{n-1}_{+}}\left(p-1\right) (\D _{\theta} w)^2 +\left(pJ_1-\frac{n(n-4)}{2} \right)|\nabla _{\theta}w|^2 +\left(pJ_2-\frac{n^2(n-4)^2}{16}\right)w^2 \; d\theta \leq 0.
\end{align*}

 If $\frac{n+4+2a}{n-4} < p < p_{JL2}(n,a)$,  implies that

 $$pJ_1-\frac{n(n-4)}{2}>0\; \,\,\mbox{and}\,\,\,\,  pJ_2-\frac{n^2(n-4)^2}{16}>0.$$
 The proof for the last inequality is very similar to tat of [\cite{hu}, Theorem 3.1], we leave the details for interested readers.
 So, it follows that $u\equiv 0$.\qed

\medskip
\section{Classification of stable solutions.}
\setcounter{equation}{0}
\medskip

For the case, $1<p\leq\frac{n+4+2a}{n-4}$ , we apply the integral estimates.
For the case, $\frac{n+4+2a}{n-4}<p<p_{JL2}(n,a)$, with the energy estimates and the desired monotonicity formula under the condition $\frac{n+4+2a}{n-4}<p<p_{JL2}(n,a)$, we can show that the stable solutions must be homogeneous solutions, hence by applying the classification of the homogeneous solutions (see Proposition \ref{prop1}), the solutions must be zero.

\medskip
\subsection {\textbf{Proof of Theorem \ref{th1.1}}}
\medskip

Since we assume that $u$ is a stable solution, then the integral estimate \eqref{2.13} holds with $C_{0}=0.$  We divide the proof in three parts.

\medskip
{\bf Step $1.$ Subcritical case: $1<p<\frac{n+4+2a}{n-4}$.}
\smallskip

 Applying \eqref{2.13}  and $1<p<\frac{n+4+2a}{n-4}$, we deduce that

$$\int_{B^{+}_{R}}v^{2}+ \int_{B^{+}_{R}}|x|^a|u|^{p+1} \leq C R^{n-\frac{4(p+1)+2a}{p-1}}\longrightarrow 0, \;\;\mbox{as} \; R\longrightarrow +\infty.$$

Consequently, we obtain $u\equiv0.$

\medskip
{\bf Step $2.$ Subcritical case: $p=\frac{n+4+2a}{n-4}$.} Applying again \eqref{2.13} we have
$$\int_{\R^{n}_{+}}v^{2}+|x|^a|u|^{p+1}<+\infty .$$

So, we get
$$\lim _{R\longrightarrow +\infty }\int_{A^{+}_{R}} v^{2}+|x|^a|u|^{p+1} \equiv0.$$

Now, applying Hölder inequality,  we derive
  \begin{align*}
R^{-4}\int_{A^{+}_{R}}u^{2}\leq CR^{-4}\left(\int_{A^{+}_{R}}|x|^{a}|u|^{p+1}\right)^{\frac{2}{p+1}}\left(\int_{A^{+}_{R}}|x|^{\frac{-2a}{p-1}}\right)^{\frac{p-1}{p+1}}.
\end{align*}
Therefore, from \eqref{2.12} we conclude then

\begin{align*}
\int_{B^{+}_{R}} v^{2}+|x|^a|u|^{p+1} \leq C R^{\left(n-\frac{2a}{p-1}\right)\frac{p-1}{p+1}-4}\left(\int_{A^{+}_{R}}|x|^a|u|^{p+1}\right)^{\frac{2}{p+1}} + C \int_{A^{+}_{R}}v^{2}.
\end{align*}

Under the assumptions  $p=\frac{n+4+2a}{n-4},$ tending $R\longrightarrow +\infty,$ we obtain  $u\equiv0.$

\medskip
{\bf Step $3.$ Supercritical case: $\frac{n+4+2a}{n-4}<p<p_{JL2}(n,a).$} We define blowing down sequences
$$u^{\lambda}(x)=\lambda^{\frac{4+a}{p-1}}u(\lambda x),\;\;\;v^{\lambda}(x)=\lambda^{\frac{4+a}{p-1}+2}v(\lambda x)\;,\; \;\;\forall\;\;\lambda>0.$$
$u^{\lambda}$ is also a smooth stable solution of \eqref{n}  on $\R^{n}_{+}$.
By rescaling \eqref{2.13},for all $\lambda >0$ and balls $B_{r}\subset \R^{n},$
\begin{eqnarray*}
\int_{B_{r}^{+}}(v^{\lambda})^{2}+|x|^a|u^{\lambda}|^{p+1} \leq Cr^{n-\frac{4(p+1)+2a}{p-1}}.
\end{eqnarray*}
In particular, $u^{\lambda}$ are uniformly bounded in $L^{p+1}_{loc}(\R^{n}_{+}).$  By elliptic estimates, $u^{\lambda}$ are also uniformly bounded in $W^{2,2}_{loc}(\R^{n}_{+}).$  Hence, up to a subsequence of $\lambda \rightarrow +\infty,$ we can assume that $u^{\lambda}\longrightarrow u^{\infty}$ weakly in $W^{2,2}_{loc}(\R^{n}_{+})\cap L^{p+1}_{loc}(\R^{n}_{+}).$ By compactness embedding, one has $u^{\lambda}\longrightarrow u^{\infty}$ strongly in $W^{2,2}_{loc}(\R^{n}_{+}).$
 Then for any ball $B_{R}^{+}(0)$, by interpolation between $L^q$ spaces and noting \eqref{2.13}, for any $q\in [1,\, p+1)$, as $\lambda \rightarrow +\infty$

\begin{eqnarray}\label{2.14}
||u^{\lambda}-u^{\infty}||_{L^{q}(B_{R}^{+}(0))} \leq ||u^{\lambda}-u^{\infty}||_{L^{1}(B_{R}^{+}(0))}^{\mu} ||u^{\lambda}-u^{\infty}||_{L^{p+1}(B_{R}^{+}(0))} ^{1-\mu}\longrightarrow0,
\end{eqnarray}

where  $\frac{1}{q}=\mu +\frac{1-\mu }{p+1}.$ That is, $u^{\lambda}\longrightarrow u^{\infty}$ in $L^{q}_{loc}(\R^{n}_{+})$ for any $q\in (1,p+1).$
\smallskip

For any function
 $\zeta\in C^{\infty}_{0}(\R^{n}_{+}),$
$$\int_{\R^{n}_{+}} \D u^{\infty} \D \zeta- |x|^a|u^{\infty}|^{p-1}u^{\infty}\zeta=\lim _{\lambda \longrightarrow \infty} \int_{\R^{n}_{+}} \D u^{\lambda} \D \zeta- |x|^a|u^{\lambda}|^{p-1}u^{\lambda}\zeta,$$
$$\int_{\R^{n}_{+}} ( \D \zeta)^{2}- p|x|^a|u^{\infty}|^{p-1}(\zeta)^{2}=\lim _{\lambda \longrightarrow \infty} \int_{\R^{n}_{+}} (\D \zeta)^{2}- p|x|^a|u^{\lambda}|^{p-1}(\zeta)^{2}\geq 0.$$
Thus $u^{\infty}\in W^{2,2}_{loc}(\R^{n}_{+})\cap L^{p+1}_{loc}(\R^{n}_{+})$ is a stable solution of \eqref{n}.

\medskip
Now, we can follow exactly the proof  of Lemmas 3.1--3.3 in Hu \cite{hu},  (see also  Lemmas 4.4--4.6 in D\'{a}vila et al.\cite{DavilaDupaigneWangWei}), to
obtain
\begin{lem}
\begin{itemize}
\item[ 1.] $ \lim_{\lambda\rightarrow +\infty}E(u,\lambda)<+\infty.$
\item[ 2.] $u^{\infty}$ is homogeneous.
\item[ 3.] $\lim_{r \rightarrow +\infty} E(u,r) =0$.
\end{itemize}
\end{lem}

\smallskip
Therefore, by the monotonicity formula we know that u is homogeneous, then $u\equiv0$,  by  Proposition \ref{prop 2.1}. This finishes the proof of Theorem \ref{th1.1}.\qed

\medskip

\section{\textbf{Classification of the finite Morse index solutions}}

\medskip

We proceed based on a Pohozaev-type identity, the decay estimates from the doubling
lemma \cite{PolacikQuittnerSouplet}, the monotonicity formula, the classification of the homogeneous solutions
and stable solutions we obtained before

\medskip
\subsection {\textbf{Subcritical and critical case}}
\medskip

Our approach consists in testing the equation  \eqref{n} against $\nabla u\cdot x \psi$ where $\psi \in C^2_c(\R^N),\;0 \leq \psi \leq 1$ is a cut-off functions satisfying

\begin{align}\label{cases}
 \left\{\begin{array}{lll}
  \psi\equiv 1\;\quad \mbox{if }\quad |x|<R,\;\;\; \quad  \psi\equiv 1\; \quad \mbox{if }\quad  |x|> 2 R,\\\\
|\nabla^q \psi|\leq C R^{-q}, \;  \quad \mbox{if }\;\;\;x \in A_R=\{R <|x|<2R\} \;\;\; q\leq 2.
\end{array}
\right.
\end{align}

 We provide the following variant of the Pohozaev identity.  In view of the cut-off functions $\psi$, we can ovoid the spherical integrals raised in  \cite{PS, WX},  which are very difficult to control. Precisely, we have
 \begin{lem}\label{lem1h.1}
Let $u$ be a solution of \eqref{n} and set $v=\Delta u$. Then for any $\psi \in C_c^{2}(\Omega)$,
 \begin{align}\label{new12}
\begin{split}
&\; \frac{n+a}{p+1}\int_{\Omega} |x|^a|u|^{p+1}\psi
- \frac{n-4}{2}\int_{\Omega} v^{2}\psi  \\
 =& \;-\frac{1}{p+1}\int_{\Omega} |x|^a|u|^{p+1}(\nabla \psi\cdot x)+ \frac{1}{2}\int_{\Omega} (\nabla\psi\cdot x)v^{2}\\
& \; - \int_{\Omega} \Big[2v(\nabla u\cdot \nabla \psi)  + 2v \nabla^2u(x, \nabla\psi) + v(\nabla u\cdot x)\Delta \psi \Big].
\end{split}
 \end{align}
 \end{lem}
\noindent{\bf Proof.} Let $\psi\in C_c^{2}(\Omega)$, multiplying equation \eqref{n} by $\nabla u\cdot x \psi$ and integrating by parts, we get
\begin{align*}
& \int_{\Omega}|x|^a|u|^{p-1}u (\nabla u\cdot x)\psi \\ = & \;  \int_{\Omega}\Delta u\D(\nabla u\cdot x\psi)
=  \; \int_{\Omega}\;v \Big[(\nabla(v)\cdot x)\psi+2v \psi+2\nabla(\nabla u\cdot x)\cdot\nabla \psi+(\nabla u\cdot x)\Delta\psi \Big].
\end{align*}
Direct calculation yields $\nabla(\nabla u\cdot x)\cdot\nabla \psi = \nabla^{2}u (x, \nabla \psi) + (\nabla u\cdot \nabla \psi)$ and
\begin{align*}
 \int_{\Omega}\; v\Big[(\nabla(v)\cdot x)\psi+2v \psi\Big]
   &= \int_{\Omega} \frac{\nabla (v^2)}{2}\cdot x \psi+ 2\int_{\Omega}v^{2}\psi \\
   &=\frac{4-n }{2}\int_{ \Omega}v^2\psi -\frac{1}{2}\int_{\Omega}v^2(\nabla \psi\cdot x).
\end{align*}
Moreover,
\begin{align*}
\int_\Omega |x|^a|u|^{p-1}u (\nabla u\cdot x)\psi
=-\frac{n+a}{p+1} \int_{\Omega}|u|^{p+1} \psi - \frac{1}{p+1} \int_{\Omega}|x|^a|u|^{p+1} x\cdot \nabla \psi .
\end{align*}
Therefore, the claim follows by regrouping the above equalities. \qed
\medskip

We claim then

\begin{lem}\label{3.3}
Let $u\in C^4(\overline{\R^n_+})$ be a solution of \eqref{n} which is stable a compact set of $\R^n_+.$ If $p\in (1,\frac{n+4+2a}{n-4}),$
then $|x|^{\frac{a}{p+1}}u \in L^{p+1}(\R^n_+),$ $v \in L^2(\R^n_+),$
\begin{align}\label{t1}
\frac{n-4}{2}\int_{\R^{n}_{+}}v^{2}=\frac{n+a}{p+1}\int_{\R^{n}_{+}}|x|^a|u|^{p+1}.
\end{align}
and
\begin{align}\label{t2}
\int_{\R^{n}_{+}}v^{2}=\int_{\R^{n}_{+}}|x|^a|u|^{p+1}.
\end{align}
\end{lem}
\noindent{\bf Proof.} Using \eqref{2.13} and tending $R \to \infty$, we obtain
\begin{align}
\label{new11}
|x|^{\frac{a}{p+1}}u \in L^{p+1}(\R^n_+) \quad \mbox{and}\quad v \in L^{2}(\R^n_+).
\end{align}

By Hölder’s inequality, there holds
$$R^{-4}\int_{A^{+}_{R}}|u|^{2} \leq C R^{(n-\frac{4(p+1)+2a}{p-1})\frac{p-1}{p+1}}\left( \int_{A^{+}_{R}}|x|^a|u|^{p+1}\right)^{\frac{2}{p+1}}.$$
On the other hand, by standard scaling argument, there exists $C > 0$ such that for any $R > 0$, any $u \in C^{4}(A^{+}_{R})$ with $A^{+}_{R}= B^{+}_{2R}\backslash B^{+}_{R}$,
\begin{align*}
R^{-2} \int_{A^{+}_{R}} |\nabla u|^2  \leq C \int_{A^{+}_{R}} v^{2}  + CR^{-4} \int_{A^{+}_{R}} u^2.
\end{align*}
Therefore, as $p$ is subcritical, we deduce that
\begin{align}
\label{new9}
 CR^{-4} \int_{A^{+}_{R}} u^2 + R^{-2} \int_{A^{+}_{R}} |\nabla u|^2 \to 0 \quad \mbox{as }\; R \to \infty.
\end{align}

Now we shall estimate the integral
$$\int_{A^{+}_{R}}| \nabla^{2} u|^{2}.$$

Since  $u\zeta=0$ on $\partial\R^n_+,$ by standard elliptic theory, there exists $C > 0$  such that
 \begin{align}\label{0.255}
 \begin{split}
\int_{A^{+}_{R}} |\nabla^2(u\zeta)|^2   \leq
C \int_{A^{+}_{R}}|\D (u \zeta)|^2  \leq C\int_{A^{+}_{R}} \Big[u^2 |\D\zeta|^2 + |\nabla u|^2|\nabla \zeta|^2 + v^2\Big] .
\end{split}
\end{align}
So, we get

 \begin{align}\label{new4}
 \begin{split}
\int_{A^{+}_{R}}| \nabla^{2} u|^{2}\zeta^{2} 
 \leq& \;C\int_{A^{+}_{R}}|\nabla^2(u\zeta)|^2 +C\int_{A^{+}_{R}}|\nabla u|^{2}|\nabla \zeta|^{2}+ C\int_{A^{+}_{R}}u^{2}\left(|\nabla \zeta|^{4} +|\nabla^{2} \zeta|^{2}\right)\\
 \leq & \; C\int_{A^{+}_{R}}v^{2} + CR^{-4} \int_{A^{+}_{R}} u^2 + R^{-2} \int_{A^{+}_{R}} |\nabla u|^2.
 \end{split}
\end{align}
 Using \eqref{new11}--\eqref{new9}, there holds
\begin{align}
\label{new10}
\int_{\R^n_+}|\nabla^2 u|^2  < \infty.
\end{align}

Now, to prove \eqref{t1}, we will show that any terms on the right hand side of \eqref{new12} (Denote by $I_R$),  tends to $0$ as $R\longrightarrow +\infty$. Remark that $\nabla \psi \ne 0$ only in $A^{+}_{R} = B^{+}_{2R}\backslash B^{+}_R$ and $\|\nabla^k \psi\|_\infty \leq C_k R^{-k}$, there holds
\begin{align*}
|I_R| \leq C\int_{A^{+}_{R}} \Big(|x|^a|u|^{p+1} + v^2\Big) + \frac{C}{R}\int_{A^{+}_{R}}|v||\nabla u| + C\int_{A^{+}_{R}}|v||\nabla^2 u|
\end{align*}
Thanks to the estimates \eqref{new11}-\eqref{new10} and H\"older's inequality, clearly $\lim_{R\to \infty}I_R = 0$, hence we get \eqref{t1}.

\medskip
On the other hand, using $u \psi$ as test function in \eqref{n}, we have
\begin{align*}
 \int_{B^{+}_{2R}}v^2 \psi-\int_{B^{+}_{2R}} |x|^a|u|^{p+1} \psi\leq C\int_{B^{+}_{2R}}|uv|
  |\D \psi|
  + C\int_{B^{+}_{2R}}|v||\nabla u||\nabla\psi| dx\leq \frac{C}{R^2}\int_{A^{+}_R}|uv|
  + \frac{C}{R}\int_{A^{+}_R}|v||\nabla u|.
\end{align*}
Apply H\"older's inequality, \eqref{new11}--\eqref{new9} and tending $R$ to $\infty$, so we obtain \eqref{t2}. The proof is completed. \qed

\bigskip
\noindent
{\bf Proof of Theorem \ref{th1.2}}.

\medskip

\medskip
{\bf Step $1.$ Subcritical case: $1<p<\frac{n+4+2a}{n-4}$.}
\smallskip

 Combining \eqref{t1} and \eqref{t2}, there holds
\begin{align*}
\left(\frac{n-4 }{2}-\frac{n+a}{p+1}\right)\int_{A^{+}_R}|u|^{p+1}=0.
\end{align*}
We are done, since $n<\frac{4(p+1)+2a}{p-1}$ implies that $\frac{n-4 }{2}-\frac{n+a}{p+1} < 0$.
\medskip

{\bf Step $2.$ Subcritical case: $p=\frac{n+4+2a}{n-4}$.}
\smallskip

We can proceed as in the proof of equality \eqref{t2}, to derive that $$\int_{\R^n_+} v^2 =\int_{\R^n_+} |x|^a|u|^{p+1} < + \infty.$$ \qed

\medskip
\subsection {\textbf{Supercritical case}}
\medskip

To classify finite Morse index solutions in the supercritical case, applying the doubling
lemma in \cite{PolacikQuittnerSouplet}, we get the following crucial lemma.

\begin{lem}\label{4.l2}
Let $n\geq 1,\;\; 1 < p < p_{JL2}(n,0)$ and $\tau\in (0, 1].$ Let $c\in C^{\tau}(\overline{B^+_{1}})$ satisfy
\begin{eqnarray}\label{4.1}
\Vert c\Vert_{C^{\tau}(\overline{B^+_{1}})}\leq C_{1}\;\;\;\mbox{and}\;\;c(x)\geq C_{2},\;\;x\in \overline{B^+_{1}},
\end{eqnarray}
for some constants $C_1, C_2 > 0.$ There exists a constant $C,$ depending only on $\alpha, C_1, C_2, p, n,$ such that, for any stable solution $u$ of
\begin{eqnarray}\label{4.2}
\D^{2}u=c(x)|u|^{p-1}u\;\;\mbox{in}\;B^+_{1} \;\mbox{and} \; u =\frac{\partial u}{\partial x_n}=0\;\;\mbox{on}\;\partial B^+_{1},
\end{eqnarray}
u satisfies $$|u(x)|^{\frac{p-1}{4}} \leq C(1+dist^{-1}(x,\partial B^+_{1})).$$
\end{lem}
\noindent{\bf Proof.} Arguing by contradiction, we suppose that there exist sequences $c_k, u_k$ verifying \eqref{4.1}-\eqref{4.2} and points $y_k,$ such that the functions $$M_k=|u_k|^{\frac{p-1}{4}}$$
satisfy $$M_k(y_k)>2k(1+dist^{-1}(y_k,\partial B^+_{1}))\geq 2kdist^{-1}(y_k,\partial B^+_{1})).$$
By the doubling lemma in \cite{PolacikQuittnerSouplet}, there exists $x_k$ such that
$$M_k(x_k)\geq M_k(y_k),\;\;\;\;M_k(x_k)\geq 2kdist^{-1}(x_k,\partial B^+_{1})),$$
and
\begin{eqnarray}\label{4.3}
M_{k}(z)\leq 2M_{k}(x_{k}), \;\;\mbox{for all}\; z\in B^+_{1} \;\mbox{such that}\;|z-x_{k}|\leq k M^{-1}_{k} (x_{k}).
\end{eqnarray}
We have
\begin{eqnarray}\label{4.4}
\lambda_k=M^{-1}_k (x_k)\longrightarrow 0,\;\;k\longrightarrow \infty
\end{eqnarray}
due to $M_k(x_k)\geq M_k(y_k)>2k.$
\smallskip

Next we let $$v_k(y)=\lambda_{k}^{\frac{4}{p-1}}u_k(x_k+\lambda_{k}y)\;\;\mbox{and}\;\;\tilde{c}_{k}(y)=c_{k}(x_k+\lambda_{k}y),\; for\; y\in B_k\;\mbox{and}\;y_n>-\frac{y_{k,n}}{\lambda_{k}},$$
where $y_k=(y_{k,1},...,y_{k,n}).$ Then, $v_k(y)$ is the solution of
\begin{equation}
\begin{cases} \D^{2} v_k(y)= \tilde{c}_{k}(y)|v_k(y)|^{p-1}v_k(y) \; ,\;\;\;|y|<k \; ,\;\;\;y_n>-\frac{y_{k,n}}{\lambda_{k}}, \\\\ v_k(y) =\frac{\partial v_k(y)}{\partial y_n} = 0 \;,\;\;\;  |y|<k \; ,\;\;\;y_n=-\frac{y_{k,n}}{\lambda_{k}} ,  \end{cases}
\end{equation}
with $$|v_k(0)|=1\;\;\mbox{and}\;\;|v_k(y)| \leq 2^{\frac{4}{p-1}}\; ,\;\;\;|y|<k \; ,\;\;\;y_n>-\frac{y_{k,n}}{\lambda_{k}}.$$
Two cases may occur as $k\longrightarrow \infty,$ either case (1)  $$\frac{y_{k,n}}{\lambda_{k}}\longrightarrow +\infty$$ for a subsequence still denoted as before, or case (2)  $$\frac{y_{k,n}}{\lambda_{k}}\longrightarrow c\geq 0.$$
If case (1), after extracting a subsequence, $\tilde{c}_{k}\longrightarrow C$ in $C_{loc}(\R^n)$ with $C>0$ a constant and we may also assume that $v_k\longrightarrow v$ in $C^4_{loc}(\R^n)$, and $v$ is a stable solution of $$\D^2 v=C|v|^{p-1}v\; \mbox{in} \;\R^n \;\;\mbox{and}\;|v(0)|=1$$
By the Liouville type Theorems in \cite{DavilaDupaigneWangWei} for stable solutions, we derive that $v\equiv 0$. This is a contradiction.
\smallskip

If case (2) we can prove that $c>0,$ thus we get a stable solution of \eqref{n} in $\overline{\R^{n}_{+}}$ and $|v(c)|=1$, which contradict Theorem \ref{th1.1} for $1 < p < p_{0}(n,4)$
\medskip

\begin{prop}\label{l5}
Let $u$ be a (positive or sign changing) solution to \eqref{n} which is stable outside a compact set of $\R^n_+$. There exist constants $C$ and $R_0$ such that
\begin{eqnarray}\label{5.8}
|u(x)| \leq C |x|^{-\frac{4+a}{p-1}},\;\;\;\mbox{for all}\; x \in B^+_{R_0}(0)^c,
\end{eqnarray}
\begin{eqnarray}\label{e036}
\sum_{k\leq 3} |x|^{\frac{4+a}{p-1}+k} |\nabla^k u(x)|\leq C ,\;\;\;\mbox{for all}\; x\in B^+_{3R_0}(0)^c .
\end{eqnarray}
\end{prop}
\noindent{\bf Proof.} Assume that $u$ is stable outside $B^+_{R_0}$ and $|x_0|>2R_0.$ We denote
$$R=\frac{1}{2}|x_0|$$
and observe that, for all $y\in B^+_1,\;\frac{|x_0|}{2}<|x_0+Ry|<\frac{3|x_0|}{2},$ so that $x_0+Ry\in B^+_{R_0}(0)^c$. Let us thus define
$$U(y)=R^{\frac{4+a}{p-1}} u(x_0+Ry).$$ Then $U$ is a solution of
$$\D^2 U=c(y)|U|^{p-1}U \;\;\mbox{in}\; B^+_1 \;\;\mbox{and}\;U=\frac{\partial U}{\partial y_n}=0\;\mbox{on}\;\partial B^+_1,\;\;\mbox{with}\;\; c(y)=\left| y+\frac{x_0}{R}\right|^{a}.$$
Notice that $| y+\frac{x_0}{R}|\in [1,3]$ for all $y\in \overline{B^+_1}.$ Moreover $\Vert c\Vert_{C^{1}(\overline{B^+_1})}\leq C(a).$ Then applying Lemma \ref{4.l2}, we have $|U(0)|\leq C ,$ hence
$$|u(x_0)|\leq CR^{-\frac{4+a}{p-1}},$$
which yields the inequality \eqref{5.8}.
\smallskip

Next,we only prove the inequality \eqref{e036}.For any $x_0$ with $|x_0|> 3 R_0$, take $\lambda= \frac{|x_0|}{2}$ and define
$$\overline{u}(x)= \lambda^{\frac{4+a}{p-1}}u(x_0+ \lambda x).$$
From \eqref{5.8}, $|\overline{u}|\leq C_0$ in $B^+_1(0)$. Standard elliptic estimates give $$\sum_{k\leq 5}|\nabla^k \overline{u}(0)|\leq C.$$
Rescaling back we get \eqref{e036}.\qed

\bigskip
\noindent
{\bf Proof of Theorem \ref{th1.2}}. {\bf Supercritical case: $p>\frac{n+4+2a}{n-4}$  and $p<p_{JL2}(n,0)$.}

\smallskip

\begin{lem}
There exists a constant $C_2,$ such that forall $r>3R_0,\;\;E(u,r)\leq C_2.$
\end{lem}
\noindent{\bf Proof.} From the monotonicity formula, combining the derivative estimates \eqref{e036}, we have then
\begin{align*}
 & \;E(u,r)\\
\leq & \; C r^{\frac{4(p+1)+2a}{p-1}-n}\left( \int_{B^{+}_r} v^2 + |x|^a|u|^{p+1} \right) +C r^{\frac{8+2a}{p-1}+1-n}\int_{ \partial B^{+}_{r}}  u^2 + C r^{\frac{8+2a}{p-1}+2-n}\int_{\partial B^{+}_{r}}|u||\nabla u| \\
 +& \; C r^{\frac{8+2a}{p-1}+3-n}\int_{\partial B^{+}_{r}} |\nabla u|^2  + C r^{\frac{8+2a}{p-1}+3-n}\int_{\partial B^{+}_{r}} |u||\nabla^2 u| + C r^{\frac{8+2a}{p-1}+4-n}\int_{\partial B^{+}_{r}} |u||\nabla^2 u| \leq  C.
 \end{align*}
 where $C$ depends on the constant appeared in \eqref{e036}. \qed

\medskip
We claim then

\begin{cor}\label{cor2}
 $$\int_{(B^{+}_{3R_0}(0))^c} \frac{\left(\frac{4+a}{p-1} |x|^{-1} u(x) + \frac{\partial u}{\partial r}(x)\right)^2}{|x|^{n-2-\frac{8+2a}{p-1}}} < + \infty.$$
 \end{cor}
\smallskip

As before, we define a blowing down sequence

$$u^{\lambda}(x)= \lambda^{\frac{4+a}{p-1}} u(\lambda x).$$

 By Proposition \ref{l5}, $u^{\lambda}$ are uniformly bounded in $C^5(B^{+}_r(0)\backslash B^{+}_{1/r}(0))$ for any fixed $r>1$. 
 
 $u^{\lambda}$ is stable outside $B^{+}_{R_0/\lambda}(0)$. There exists a function $u^{\infty}\in C^6(\R^n\backslash \{0 \})$, such that up to a subsequence of $\lambda \rightarrow +\infty$, $u^{\lambda}$ converges to $u^{\infty}\in C_{loc}^4(\R^n_{+}\backslash \{0 \})$. $u^{\infty}$ is a stable solution of \eqref{n} in $\R^n_{+}\backslash \{0 \}$.
\smallskip

 Using  Corollary \ref{cor2}, we obtain for any $r>1$, 
 \begin{align*}
 & \;\int_{B^{+}_r\backslash B^{+}_{1/r}} \frac{\left(\frac{4+a}{p-1} |x|^{-1} u^{\infty}(x) + \frac{\partial u^{\infty}}{\partial r}(x)\right)^2}{|x|^{n-2-\frac{8+2a}{p-1}}} \\=& \; \lim_{\lambda\rightarrow +\infty} \int_{B^{+}_r\backslash B^{+}_{1/r}} \frac{\left(\frac{4+a}{p-1} |x|^{-1} u^{\lambda}(x) + \frac{\partial u^{\lambda}}{\partial r}(x)\right)^2}{|x|^{n-2-\frac{8+2a}{p-1}}} 
 = \;\lim_{\lambda\rightarrow +\infty} \int_{B^{+}_r\backslash B^{+}_{1/r}} \frac{\left(\frac{4+a}{p-1} |x|^{-1} u(x) + \frac{\partial u}{\partial r}(x)\right)^2}{|x|^{n-2-\frac{8+2a}{p-1}}}=0.
  \end{align*}
  
  Hence, $u^{\infty}$ is homogeneous, and from Proposition \ref{prop1}, $u^{\infty}\equiv0$. This holds for every limit of $u^{\lambda}$ as $\lambda\rightarrow +\infty$, thus we get  
  
  $$\lim_{|x|\rightarrow +\infty} |x|^{\frac{4+a}{p-1}}|u(x)|=0.$$
  
From \eqref{e036}, we derive 

$$\lim_{|x|\rightarrow+\infty} \sum_{k\leq 4}|x|^{\frac{4+a}{p-1}+k} |\nabla^k u(x)|=0.$$

For $\varepsilon > 0$, take an $R$ such that for $|x|> R$,
$$\sum_{k\leq 4}|x|^{\frac{4+a}{p-1}+k} |\nabla^k u(x)|\leq \varepsilon.$$

Then for $r>>R$,
\begin{align*}
E(u,r)&\leq C r^{\frac{4(p+1)+2a}{p-1}-n}\left( \int_{B^{+}_R(0)} v^2 + |x|^a|u|^{p+1}\right) + C \epsilon r^{\frac{8+2a}{p-1}+4-n}\int_{B^{+}_r(0) \backslash B^{+}_R(0)} |x|^{-\frac{8+2a}{p-1}-4}\\& \;
+ C \epsilon r^{\frac{8+2a}{p-1}+5-n} \int_{\partial B^{+}_r(0)} |x|^{-\frac{8+2a}{p-1}-4}\;\; \leq C(R) \left( r^{\frac{4(p+1)+2a}{p-1}-n}+ \varepsilon\right).
\end{align*}
Since $ \frac{4(p+1)+2a}{p-1}-n < 0$ and $\varepsilon$ can be arbitrarily small, we derive $\lim_{r\rightarrow +\infty}E(u,r)=0$. Because $\lim_{r \rightarrow 0} E(r,u)=0$ ( by the smoothness of $u$ ), the same argument for stable solutions implies that $u\equiv 0$. \qed

\medskip
\section*{Acknowledgment}

The authors extend their appreciation to the Deanship of Scientific Research at King Khalid University, Abha, KSA for funding this work through Research Group under grant number (R.G.P-2 / 121/ 42).

\end{document}